\documentclass[thmsa]{article}
\usepackage{amssymb}

\input{tcilatex}
\input{diagrams.tex}
\newarrow{Onto}{-}{-}{-}{-}{->>}
\newarrow{Into}{C}{-}{-}{-}{>}
\newarrow{Corr}{<}{-}{-}{-}{>}
\newarrow{Ev}{|}{-}{-}{-}{>}

\begin{document}

\author{Sergio D. Grillo \\
\textit{Centro At\'{o}mico Bariloche and Instituto Balseiro}\\
\textit{\ 8400-S. C. de Bariloche}\\
\textit{\ Argentina}}
\title{\emph{FRT} Construction and Equipped Quantum Linear Spaces}
\date{February 2003 }
\maketitle

\begin{abstract}
We show there exists a rigid monoidal category formed out by quantum spaces
with an additional structure, such that \emph{FRT }bialgebras and
corresponding rectangular generalizations are its internal coEnd and coHom
objects, respectively. This enable us to think of them as the coordinate
rings of `quantum spaces of homomorphisms' that preserve the mentioned
structure. The well known algebra epimorphisms between \emph{FRT} bialgebras
and Manin quantum semigroups translate into `inclusions' of the
corresponding quantum spaces, as the space of endomorphisms of a metric
linear space $\mathbf{V}$ is included in $gl\left( \mathbf{V}\right) $. Our
study is mainly developed for quadratic quantum spaces, and later
generalized to the conic case.
\end{abstract}

\section{Introduction}

Given a finite dimensional $\Bbbk $-vector space $\mathbf{V}$ and a linear
endomorphism $\Bbb{R}$ of $\mathbf{V\otimes V}$, a universal bialgebra $%
A\left( \Bbb{R}\right) $ can be constructed \cite{kas}. The assignment of $%
A\left( \Bbb{R}\right) $ to each pair $\left( \mathbf{V},\Bbb{R}\right) $ is
known as \emph{FRT} \emph{construction} \cite{frt}. Every $A\left( \Bbb{R}%
\right) $ is a generically non commutative quadratic algebra generated by a
finite dimensional coalgebra, more precisely, by linearly independent
coefficients of a multiplicative matrix $\mathbf{t}$ \cite{A}. This is why
they are called \emph{quantum matrix bialgebras}. Given a basis $\left\{
v_{i}\right\} $ of $\mathbf{V}$, coefficients of $\mathbf{t}$ are identified
with the elements $t_{i}^{j}=v^{j}\otimes v_{i}$ of $\mathbf{V}^{\ast
}\otimes \mathbf{V}$, and the quadratic relations they must satisfy are
sometimes written 
\[
\Bbb{R}_{ij}^{kl}\;t_{k}^{n}\,t_{l}^{m}-t_{i}^{k}\,t_{j}^{l}\;\Bbb{R}%
_{kl}^{nm};\;\;\;i,j,n,m=1,...,\dim \mathbf{V,} 
\]
being $\Bbb{R}_{ij}^{kl}\in \Bbbk $ the coefficients of $\Bbb{R}$ in the
given basis.

Quantum matrix bialgebras are the `dual' version of quantum universal
enveloping algebras, such as Drinfeld-Jimbo \cite{dj} quantized Lie
bialgebras $U_{q}\left( \frak{g}\right) $.\footnote{%
It is worth mentioning we are not asking for $\Bbb{R}$ to be a Yang-Baxter
operator. If this were the case, $\Bbb{R}$ would indicate the so called $R$%
-matrix.} On the other hand, they are quotient of Manin \emph{quantum }(%
\emph{semi})\emph{groups} $\underline{end}\left[ \mathcal{V}\right] $ \cite
{man0}, i.e. the internal coEnd objects of the monoidal category $\mathrm{QA}
$ of quadratic algebras. In other words, there exists a bialgebra
epimorphism $\underline{end}\left[ \mathcal{V}\right] \twoheadrightarrow
A\left( \Bbb{R}\right) $ in $\mathrm{QA}$. The relationship between $\Bbb{R}$
and $\mathcal{V}$ will be discussed later. By now, let us write $\mathcal{V}%
\Vdash \Bbb{R}$ when they are related. Using geometric language, each object 
$\underline{end}\left[ \mathcal{V}\right] $ is interpreted as the coordinate
ring of a non commutative algebraic variety, or \emph{quantum linear space},
living in the opposite category $\mathrm{QA}^{op}$. It represents the
quantum semigroup of endomorphisms corresponding to the quantum space $%
\mathcal{V}^{op}$. Thus, epimorphism above gives rise to a monic $A\left( 
\Bbb{R}\right) ^{op}\hookrightarrow \underline{end}\left[ \mathcal{V}\right]
^{op}$ enabling us to regard $A\left( \Bbb{R}\right) $ as the coordinate
ring of a quantum subspace of $\underline{end}\left[ \mathcal{V}\right]
^{op} $. Of course, Manin construction also includes \emph{quantum spaces of
homomorphisms} $\underline{hom}\left[ \mathcal{W},\mathcal{V}\right] $,
indicating by $\mathcal{W}$ a quadratic algebra generated by a vector
subspace $\mathbf{W}$. They have as \emph{FRT} analogue the \emph{%
rectangular quantum matrix algebras} $A\left( \Bbb{R}:\Bbb{S}\right) $ \cite
{mm}\cite{maj}, in the sense that there exist algebra epimorphisms $%
\underline{hom}\left[ \mathcal{W},\mathcal{V}\right] \twoheadrightarrow
A\left( \Bbb{R}:\Bbb{S}\right) $ leading us to a geometric interpretation as
described before. $\Bbb{S}$ denotes a linear map $\mathbf{W}\otimes \mathbf{W%
}\rightarrow \mathbf{W}\otimes \mathbf{W}$ such that $\mathcal{W}\Vdash \Bbb{%
S}$. Given a basis $\left\{ w_{i}\right\} $ of $\mathbf{W}$, the algebras $%
A\left( \Bbb{R}:\Bbb{S}\right) $ are generated by symbols $%
t_{i}^{j}=w^{j}\otimes v_{i}\in \mathbf{W}^{\ast }\otimes \mathbf{V}$
satisfying\footnote{%
We are using the convention of ref. \cite{phh} to evaluate the maps $\Bbb{R}$
and $\Bbb{S}$ appearing in the rectangular quantum matrix algebras, instead
of the one used in \cite{mm}.} 
\begin{equation}
\Bbb{R}_{ij}^{kl}\;t_{k}^{n}\,t_{l}^{m}-t_{i}^{k}\,t_{j}^{l}\;\Bbb{S}%
_{kl}^{nm};\;\;\;i,j=1,...,\dim \mathbf{V},\;\;\;n,m=1,...,\dim \mathbf{W}.
\label{rqa}
\end{equation}
These algebras were studied in detail in \cite{phh}, where $\Bbb{R}$ and $%
\Bbb{S}$ are Yang-Baxter operators of Hecke type.

This paper was mainly motived by the following question induced by
`inclusions' $A\left( \Bbb{R}:\Bbb{S}\right) ^{op}\hookrightarrow \underline{%
hom}\left[ \mathcal{W},\mathcal{V}\right] ^{op}$: Do the quadratic algebras $%
A\left( \Bbb{R}:\Bbb{S}\right) $ represent homomorphisms between quantum
spaces supplied with some additional structure, i.e. spaces that are not
characterized just by its respective coordinate rings ?

In order to answer this question we encode Manin and \emph{FRT}
constructions, reformulating and generalizing the latter, in the unifying
language of rigid monoidal categories \cite{dm}. We show that bialgebras $%
A\left( \Bbb{R}\right) $ can be seen as internal coEnd objects contained in
certain rigid monoidal category $\left( \mathrm{EQA},\boxtimes \right) $,
the \emph{equipped quantum spaces}, formed out by pairs $\frak{V}=\left( 
\mathcal{V};\Bbb{R}\right) $ with $\mathcal{V}\in \mathrm{QA}$ and $\mathcal{%
V}\Vdash \Bbb{R}$. More precisely, there exists a surjective embedding $%
\mathsf{U}:\mathrm{EQA}\hookrightarrow \mathrm{QA}:\left( \mathcal{V};\Bbb{R}%
\right) \mapsto \mathcal{V}$, and a related opposite $\mathsf{U}^{op}:%
\mathrm{EQA}^{op}\hookrightarrow \mathrm{QA}^{op}$, such that for each pair $%
\frak{V}$, the object $\underline{hom}\left[ \frak{V},\frak{V}\right] =%
\underline{end}\left[ \frak{V}\right] $ is \emph{functored} to $A\left( \Bbb{%
R}\right) $. In general, coHom objects $\underline{hom}\left[ \frak{W},\frak{%
V}\right] $, with $\frak{W}=\left( \mathcal{W};\Bbb{S}\right) $, are sent to 
$A\left( \Bbb{R}:\Bbb{S}\right) $. Moreover, from the general formalism of
rigid monoidal categories follows existence and associativity properties of 
\emph{rectangular comultiplication} maps defined in \cite{mm}, and also
existence of algebra epimorphisms 
\begin{equation}
\underline{hom}_{\mathrm{QA}}\left[ \mathsf{U}\frak{W},\mathsf{U}\frak{V}%
\right] \twoheadrightarrow \mathsf{U}\underline{hom}_{\mathrm{EQA}}\left[ 
\frak{W},\frak{V}\right] ,\;\;\;\forall \frak{W},\frak{V}\in \mathrm{EQA},
\label{inc}
\end{equation}
from which previously mentioned `inclusions' are deduced. We conclude, each
algebra $A\left( \Bbb{R}:\Bbb{S}\right) $ is the coordinate ring of the
space $\underline{hom}\left[ \frak{W},\frak{V}\right] ^{op}\in \mathrm{EQA}%
^{op}$ of homomorphisms between spaces $\frak{W}^{op}$ and $\frak{V}^{op}$.
Thus, such spaces are described by their respective coordinate rings $%
\mathsf{U}\underline{hom}\left[ \frak{W},\frak{V}\right] =A\left( \Bbb{R}:%
\Bbb{S}\right) $, $\mathsf{U}\frak{W}$ and $\mathsf{U}\frak{V}$ (given by
quadratic algebras), and an additional data. We also show $\mathrm{EQA}$ is
equivalent to a category whose objects are pairs $\left( \mathbf{V},\Bbb{R}%
\right) $ as described above, in such a way that we can write $\left( 
\mathcal{V};\Bbb{R}\right) \equiv \left( \mathbf{V},\Bbb{R}\right) $ if $%
\mathcal{V}$ is generated by $\mathbf{V}$. Hence, we are assigning a
bialgebra $A\left( \Bbb{R}\right) =\mathsf{U}\underline{end}\left[ \frak{V}%
\right] $ to each pair $\left( \mathbf{V},\Bbb{R}\right) $ in a universal
way.

\section{Quantum linear spaces}

In what follows $\Bbbk $ indicates some of the numerics fields, $\Bbb{R}$ or 
$\Bbb{C}$. The usual tensor product on $\Bbbk \mathrm{-Alg}=\mathrm{Alg}$%
\textrm{\ }and $\mathrm{Vct}_{\Bbbk }=\mathrm{Vct}$ (the categories of
unital associative $\Bbbk $-algebras and of $\Bbbk $-vector spaces,
respectively) is denoted by $\otimes $. \textrm{Vct}$_{f}$ indicates the
full subcategory of $\mathrm{Vct}$ formed out by finite dimensional vector
spaces.

Originally \cite{man0}, Manin defined quantum spaces as opposite objects to
quadratics algebras. The latter are pairs $\left( \mathbf{A}_{1},\mathbf{A}%
\right) $, with $\mathbf{A}\in \mathrm{Alg}$ generated by $\mathbf{A}_{1}$
in $\mathrm{Vct}_{f}$, such that the canonical epimorphism $\mathbf{A}%
_{1}^{\otimes }\twoheadrightarrow \mathbf{A}$ has as kernel a bilateral
ideal algebraically generated by a subspace of $\mathbf{A}_{1}^{\otimes 2}$.
As usual, $\mathbf{A}_{1}^{\otimes }=\tbigoplus_{n\in \Bbb{N}_{0}}\mathbf{A}%
_{1}^{\otimes n}$ denotes the tensor algebra of\textbf{\ }$\mathbf{A}_{1}$
(being $\Bbb{N}_{0}\doteq \Bbb{N}\cup \left\{ 0\right\} $). To be more
explicit, for every quadratic algebra $\left( \mathbf{A}_{1},\mathbf{A}%
\right) $ there exists a subspace $\mathbf{R}\subset \mathbf{A}_{1}^{\otimes
2}$ such that 
\[
\ker \left[ \mathbf{A}_{1}^{\otimes }\twoheadrightarrow \mathbf{A}\right] =I%
\left[ \mathbf{R}\right] =\mathbf{A}_{1}^{\otimes }\cdot \mathbf{R}\cdot 
\mathbf{A}_{1}^{\otimes }. 
\]
In general, we note by $I\left[ \mathbf{X}\right] \subset \mathbf{A}%
_{1}^{\otimes }$ the bilateral ideal generated by a set $\mathbf{X}\subset 
\mathbf{A}_{1}^{\otimes }$. For instance, each algebra $A\left( \Bbb{R}:\Bbb{%
S}\right) $ defines a quantum space 
\[
A\left( \Bbb{R}:\Bbb{S}\right) \equiv \left( \mathbf{W}^{\ast }\otimes 
\mathbf{V},A\left( \Bbb{R}:\Bbb{S}\right) \right) . 
\]
The kernel of its related canonical epimorphism is generated by the elements
given in Eq. $\left( \ref{rqa}\right) $. The category $\mathrm{QA}$, as
mentioned before, has above pairs as objects and as arrows $\left( \mathbf{A}%
_{1},\mathbf{A}\right) \rightarrow \left( \mathbf{B}_{1},\mathbf{B}\right) $
algebra homomorphisms $\mathbf{A}\rightarrow \mathbf{B}$ that preserve
generating spaces, that is to say, $\mathbf{A}\rightarrow \mathbf{B}$
restricted to $\mathbf{A}_{1}$ defines a linear map $\mathbf{A}%
_{1}\rightarrow \mathbf{B}_{1}$.

In \cite{man1}, Manin extended the concept to arbitrary finitely generated
algebras, i.e. pairs $\left( \mathbf{A}_{1},\mathbf{A}\right) $ as above,
but without restrictions on their respective canonical epimorphisms $\mathbf{%
A}_{1}^{\otimes }\twoheadrightarrow \mathbf{A}$. We shall indicate \textrm{%
FGA }the category formed out by these pairs. Its arrows are again algebra
homomorphisms preserving the generating linear spaces. Thus, $\mathrm{QA}$
is a full subcategory of \textrm{FGA}. Note that arrows $\alpha :\left( 
\mathbf{A}_{1},\mathbf{A}\right) \rightarrow \left( \mathbf{B}_{1},\mathbf{B}%
\right) $ in \textrm{FGA }are characterized by linear maps $\alpha _{1}:%
\mathbf{A}_{1}\rightarrow \mathbf{B}_{1}$ such that 
\[
\alpha _{1}^{\otimes }\left( \ker \left[ \mathbf{A}_{1}^{\otimes
}\twoheadrightarrow \mathbf{A}\right] \right) \subset \ker \left[ \mathbf{B}%
_{1}^{\otimes }\twoheadrightarrow \mathbf{B}\right] .
\]
In $\mathrm{QA}$, if $\mathbf{R}$ and $\mathbf{S}$ are the subspaces
generating the respective kernels, last condition reads $\alpha
_{1}^{\otimes 2}\left( \mathbf{R}\right) \subset \mathbf{S}$. In \cite{gm}
we study another full subcategory of \textrm{FGA}, namely $\mathrm{CA}$, the
conic algebras or conic quantum spaces. Its objects $\left( \mathbf{A}_{1},%
\mathbf{A}\right) $ are such that $\mathbf{A}$ is a graded algebra and $%
\mathbf{A}_{1}$ is its subspace of homogeneous elements of degree one, or
equivalently, its related ideal (i.e. the kernel of its canonical
epimorphism) is a graded subalgebra of $\mathbf{A}_{1}^{\otimes }$. Examples
of them, beside quadratics, are the so called $m$-th quantum spaces, whose
associated ideals are generated by a subspace of $\mathbf{A}_{1}^{\otimes m}$%
, for some $m\geq 2$. The latter, in turn, form a full subcategory $\mathrm{%
CA}^{m}$ of $\mathrm{CA}$, leading us to the full inclusions $\mathrm{CA}%
^{m}\subset \mathrm{CA}\subset \mathrm{FGA}$. Of course, $\mathrm{QA}=%
\mathrm{CA}^{2}$.

The monoid we consider on these categories is the bifunctor $\circ $, given
on objects by 
\begin{equation}
\left( \mathbf{A}_{1},\mathbf{A}\right) \circ \left( \mathbf{B}_{1},\mathbf{B%
}\right) =\left( \mathbf{A}_{1}\otimes \mathbf{B}_{1},\mathbf{A}\circ 
\mathbf{B}\right) ,  \label{mon}
\end{equation}
with $\mathbf{A}\circ \mathbf{B}$ the subalgebra of $\mathbf{A}\otimes 
\mathbf{B}$ generated by $\mathbf{A}_{1}\otimes \mathbf{B}_{1}$. On arrows,
it assigns to $\alpha $ and $\beta $, with domains $\left( \mathbf{A}_{1},%
\mathbf{A}\right) $ and $\left( \mathbf{B}_{1},\mathbf{B}\right) $,
respectively, the algebra morphism $\alpha \circ \beta =\left. \alpha
\otimes \beta \right| _{\mathbf{A}\circ \mathbf{B}}$. The unit object is $%
\mathcal{I}=\left( \Bbbk ,\Bbbk \right) $ in $\mathrm{FGA}$ and $\mathcal{K}%
=\left( \Bbbk ,\Bbbk ^{\otimes }\right) $ in $\mathrm{CA}$ and every $%
\mathrm{CA}^{m}$. Let us mention the forgetful functor $\mathsf{F}:\mathrm{%
FGA}\rightarrow \mathrm{Alg}:\left( \mathbf{A}_{1},\mathbf{A}\right) \mapsto 
\mathbf{A}$ preserves the units, since $\mathsf{F}\mathcal{I}=\Bbbk $, but
is not monoidal. Nevertheless, it is easy to check that algebra inclusions $%
i_{\mathcal{A},\mathcal{B}}:\mathbf{A}\circ \mathbf{B}\hookrightarrow 
\mathbf{A}\otimes \mathbf{B}$, related to quantum spaces $\mathcal{A}=\left( 
\mathbf{A}_{1},\mathbf{A}\right) $ and $\mathcal{B}=\left( \mathbf{B}_{1},%
\mathbf{B}\right) $, define a natural transformation $\mathsf{F}\,\circ
\,\rightarrow \otimes \,\left( \mathsf{F}\times \mathsf{F}\right) $. (Note
that $\mathsf{F}\left( \mathcal{A}\circ \mathcal{B}\right) =\mathbf{A}\circ 
\mathbf{B}$ and $\mathsf{F}\mathcal{A}\otimes \mathsf{F}\mathcal{B}=\mathbf{A%
}\otimes \mathbf{B}$.) That is to say, for any couple of arrows $\alpha
,\beta \in \mathrm{FGA}$, with $\alpha :\mathcal{A}\rightarrow \mathcal{C}$
and $\beta :\mathcal{B}\rightarrow \mathcal{D}$, the diagram 
\begin{equation}
\begin{diagram}[s=2.5em] \QTR{bf}{A}\circ \QTR{bf}{B} &
\rInto^{i_{\QTR{cal}{A},\QTR{cal}{B}}} & \QTR{bf}{A}\otimes \QTR{bf}{B} \\
\dTo^{\QTR{sf}{F}\left( \alpha \circ \beta \right) }& &
\dTo_{\QTR{sf}{F}\alpha \otimes \QTR{sf}{F}\beta } \\ \QTR{bf}{C}\circ
\QTR{bf}{D} & \rInto^{i_{\QTR{cal}{C},\QTR{cal}{D}}} & \QTR{bf}{C}\otimes
\QTR{bf}{D} \\ \end{diagram}  \label{fa}
\end{equation}
is commutative. When restricted to $\mathrm{CA}$ and every $\mathrm{CA}^{m}$%
, the above natural transformation holds, but $\mathsf{F}$ does not respects
units, because $\mathsf{F}\mathcal{K}=\Bbbk ^{\otimes }$. However, the
canonical projection $\Bbbk ^{\otimes }\twoheadrightarrow \Bbbk $ defines
epimorphisms $p_{\mathcal{A}}:\Bbbk ^{\otimes }\otimes \mathbf{A}%
\twoheadrightarrow \Bbbk \otimes \mathbf{A}$, with $\mathsf{F}\mathcal{A}=%
\mathbf{A}$, that make commutative the diagrams 
\begin{equation}
\begin{diagram}[w=3.5em] \Bbbk ^{\otimes }\circ \QTR{bf}{A} &
\rInto^{i_{\QTR{cal}{K},\QTR{cal}{A}}} & \Bbbk ^{\otimes }\otimes
\QTR{bf}{A} & \rOnto^{p_{\QTR{cal}{A}}} & \Bbbk \otimes \QTR{bf}{A} \\ &
\rdTo~{\backsimeq } & & \ruTo~{\backsimeq} & \\ & & \QTR{bf}{A} & & \\
\end{diagram}  \label{fu}
\end{equation}
Indicating by $e$ the generator of $\Bbbk $, then 
\[
\Bbbk ^{\otimes }=\Bbbk \left[ e\right] \;\;\;\;and\;\;\;p_{\mathcal{A}%
}\left( e^{n}\otimes a\right) =e\otimes a.
\]
The isomorphisms $\mathbf{A}\backsimeq \Bbbk \otimes \mathbf{A}$ and $%
\mathbf{A}\backsimeq \Bbbk ^{\otimes }\circ \mathbf{A}$ are the functorial
isomorphisms related to the left unital constraint in $\mathrm{Alg}$ and $%
\mathrm{CA}$, respectively. Of course, a diagram analogous to $\left( \ref
{fu}\right) $ but with $\Bbbk $ on the right is also fulfilled.

\bigskip

There exist internal coHom objects on each one of this monoidal categories.
For instance, for $\mathcal{A}=\left( \mathbf{A}_{1},\mathbf{A}\right) $ and 
$\mathcal{B}=\left( \mathbf{B}_{1},\mathbf{B}\right) $ in $\mathrm{CA}$
(resp. $\mathrm{CA}^{m}$), they are given by graded algebras $\underline{hom}%
\left[ \mathcal{B},\mathcal{A}\right] $ generated by $\mathbf{B}_{1}^{\ast
}\otimes \mathbf{A}_{1}$ and constrained by homogeneous relations (resp. of $%
m$-th order). For more details, see \cite{gm}.

\section{The equipped quantum spaces}

Consider $\mathcal{A}=\left( \mathbf{A}_{1},\mathbf{A}\right) \in \mathrm{QA}
$ and a linear map $\Bbb{R}:\mathbf{A}_{1}^{\otimes 2}\rightarrow \mathbf{A}%
_{1}^{\otimes 2}$.

\begin{definition}
We say $\Bbb{R}$ is \textbf{compatible} with $\mathcal{A}$, and use the
shorthand notation $\mathcal{A}\Vdash $ $\Bbb{R}$, if 
\[
\ker \left[ \mathbf{A}_{1}^{\otimes }\twoheadrightarrow \mathbf{A}\right] =I%
\left[ \func{Im}\Bbb{R}\right] .
\]
A pair $\left( \mathcal{A};\Bbb{R}\right) $ such that $\mathcal{A}\Vdash $ $%
\Bbb{R}$ will be called \textbf{equipped quantum space} \emph{(}or equipped
quadratic algebra\emph{)} with \textbf{structure }$\Bbb{R}$. If a morphism
of quantum spaces $\alpha :\mathcal{A}\rightarrow \mathcal{B}$, with $%
\mathcal{A}\Vdash $ $\Bbb{R}$ and $\mathcal{B}\Vdash \Bbb{S}$, satisfies $%
\alpha _{1}^{\otimes 2}\,\Bbb{R}=\Bbb{S}\,\alpha _{1}^{\otimes 2}$, we say
that $\alpha $ preserves structures $\Bbb{R}$ and $\Bbb{S}$. The category
formed out by equipped quantum spaces and structure preserving arrows will
be denoted $\mathrm{EQA}$.\ \ \ $\blacksquare $
\end{definition}

From now on, we reserve the name \emph{pair }only for equipped quantum
spaces (in contrast to last section where we used it for ordinary ones).

A simple characterization of equipped quantum spaces is given by the
following result.

\begin{proposition}
The category $\mathrm{EQA}$ is equivalent to one whose objects are pairs $%
\left( \mathbf{V},\Bbb{R}\right) $, with $\mathbf{V}\in \mathrm{Vct}_{f}$
and $\Bbb{R}:\mathbf{V}^{\otimes 2}\rightarrow \mathbf{V}^{\otimes 2}$ a
linear map, and whose arrows $\left( \mathbf{V},\Bbb{R}\right) \rightarrow
\left( \mathbf{W},\Bbb{S}\right) $ are linear homomorphisms $l:\mathbf{V}%
\rightarrow \mathbf{W}$ such that $l^{\otimes 2}\,\Bbb{R}=\Bbb{S}%
\,l^{\otimes 2}$.
\end{proposition}

\textbf{Proof:} The equivalence is defined by functors 
\begin{equation}
\mathsf{f}:\left( \mathcal{A};\Bbb{R}\right) \mapsto \left( \mathbf{A}_{1},%
\Bbb{R}\right) \;\;\;and\;\;\;\mathsf{g}:\left( \mathbf{V},\Bbb{R}\right)
\mapsto \left( \left( \mathbf{V},\left. \mathbf{V}^{\otimes }\right/ I\left[ 
\func{Im}\Bbb{R}\right] \right) ;\Bbb{R}\right) .  \label{fg}
\end{equation}
On arrows, $\mathsf{f}\alpha =\alpha _{1}$ and $\mathsf{g}l$ is the
extension of $l$ to an algebra homomorphism. If $l$ goes from $\left( 
\mathbf{V},\Bbb{R}\right) $ to $\left( \mathbf{W},\Bbb{S}\right) $, since 
\begin{eqnarray*}
l^{\otimes }\left( I\left[ \func{Im}\Bbb{R}\right] \right) &=&\mathbf{W}%
^{\otimes }\cdot l^{\otimes }\left( \func{Im}\Bbb{R}\right) \cdot \mathbf{W}%
^{\otimes }=\mathbf{W}^{\otimes }\cdot \left( \func{Im}l^{\otimes 2}\,\Bbb{R}%
\right) \cdot \mathbf{W}^{\otimes } \\
&=&\mathbf{W}^{\otimes }\cdot \left( \func{Im}\Bbb{S}\,l^{\otimes 2}\right)
\cdot \mathbf{W}^{\otimes }\subset \mathbf{W}^{\otimes }\cdot \func{Im}\Bbb{S%
}\,\cdot \mathbf{W}^{\otimes }=I\left[ \func{Im}\Bbb{S}\right] ,
\end{eqnarray*}
then such extension is well defined.

Natural equivalence $\mathsf{f}\circ \mathsf{g}\backsimeq id$ is immediate.
The functorial isomorphisms for equivalence $\mathsf{g}\circ \mathsf{f}%
\backsimeq id$ are given by the algebra isomorphisms $\mathbf{A}\backsimeq
\left. \mathbf{A}_{1}^{\otimes }\right/ I\left[ \func{Im}\Bbb{R}\right] $,
which are well defined provided $\ker \left[ \mathbf{A}_{1}^{\otimes
}\twoheadrightarrow \mathbf{A}\right] =I\left[ \func{Im}\Bbb{R}\right]
.\;\;\;\blacksquare $

\bigskip

Because of this equivalence, we identify the objects of both categories.
That is to say, we understand pairs $\left( \mathbf{A}_{1},\Bbb{R}\right) $
and $\left( \mathcal{A};\Bbb{R}\right) $ as the same thing, indicating both
categories by \textrm{EQA}.\footnote{%
Note that the category $\mathcal{YB}$ defined in \cite{mm}, formed out by
pairs $\left( \mathbf{V},\Bbb{R}\right) $ such that $\Bbb{R}$ is a
Yang-Baxter solution of $q$-Hecke type, is a full subcategory of $\mathrm{EQA%
}$.} Since to deal with pairs $\left( \mathbf{V},\Bbb{R}\right) $ is often
easier than to deal with $\left( \mathcal{A};\Bbb{R}\right) $, we shall work
out our constructions mainly in terms of the former. Naturally, they also
provide a more direct contact to \emph{FRT }construction.

\subsection{Products and duals}

A monoidal structure and an involution can be attached to \textrm{EQA }in
the following way. Let us consider the canonical algebra isomorphism $%
\varphi _{\mathbf{V},\mathbf{W}}$ between $\left[ \mathbf{V}\otimes \mathbf{W%
}\right] ^{\otimes }$ and 
\[
\mathbf{V}\circ \mathbf{W}\doteq \bigoplus\nolimits_{n\in \Bbb{N}_{0}}%
\mathbf{V}^{\otimes n}\otimes \mathbf{W}^{\otimes n}, 
\]
the subalgebra of $\mathbf{V}^{\otimes }\otimes \mathbf{W}^{\otimes
}=\bigoplus\nolimits_{n,m}\mathbf{V}^{\otimes n}\otimes \mathbf{W}^{\otimes
m}$ generated by $\mathbf{V}\otimes \mathbf{W}$. The restriction of $\varphi
_{\mathbf{V},\mathbf{W}}$ to $\left( \mathbf{V}\otimes \mathbf{W}\right)
^{\otimes 2}$, which we also denote $\varphi _{\mathbf{V},\mathbf{W}}$, is
given by 
\[
v\otimes w\otimes v^{\prime }\otimes w^{\prime }\mapsto v\otimes v^{\prime
}\otimes w\otimes w^{\prime };\;\;\forall v,v^{\prime }\in \mathbf{V}%
,\;w,w^{\prime }\in \mathbf{W}. 
\]
We define the bifunctor $\boxtimes :\mathrm{EQA}\times \mathrm{EQA}%
\rightarrow \mathrm{EQA}$ as 
\begin{equation}
\left( \mathbf{V},\Bbb{R}\right) \times \left( \mathbf{W},\Bbb{S}\right)
\mapsto \left( \mathbf{V}\otimes \mathbf{W},\Bbb{R}\boxtimes \Bbb{S}\right)
,\;\;\;k\times l\mapsto k\boxtimes l\doteq k\otimes l,  \label{A}
\end{equation}
being 
\begin{equation}
\Bbb{R}\boxtimes \Bbb{S}\doteq \varphi _{\mathbf{V},\mathbf{W}}^{-1}\,\left( 
\Bbb{R}\otimes \Bbb{I}+\Bbb{I}\otimes \Bbb{S}\right) \,\varphi _{\mathbf{V},%
\mathbf{W}}  \label{B}
\end{equation}
(taking into account the restriction of $\varphi _{\mathbf{V},\mathbf{W}}$). 
$\Bbb{I}$ denotes the identity endomorphism of the corresponding vector
spaces. Identifying $\mathbf{V}\circ \mathbf{W}$ and $\left[ \mathbf{V}%
\otimes \mathbf{W}\right] ^{\otimes }$ we shall write $\Bbb{R}\boxtimes \Bbb{%
S}\thickapprox \Bbb{R}\otimes \Bbb{I}+\Bbb{I}\otimes \Bbb{S}$.
Straightforwardly, the bifunctor $\boxtimes $ defines a symmetric monoidal
structure with unit object $\frak{K}=\left( \Bbbk ,\Bbb{O}\right) $, being $%
\Bbb{O}$ the null endomorphism of $\Bbbk ^{\otimes 2}$, i.e. $\func{Im}\Bbb{O%
}=\left\{ 0\right\} $. The functorial isomorphisms $\tau _{\frak{V},\frak{W}%
}:\frak{V}\boxtimes \frak{W}\backsimeq \frak{W}\boxtimes \frak{V}$ related
to symmetry, with $\frak{V}=\left( \mathbf{V},\Bbb{R}\right) $ and $\frak{W}%
=\left( \mathbf{W},\Bbb{S}\right) $, are given by canonical flipping maps $%
\mathbf{V}\otimes \mathbf{W}\backsimeq \mathbf{W}\otimes \mathbf{V}$, $%
v\otimes w\mapsto w\otimes v$. The ones related to unit are $\ell _{\frak{V}%
}:v\in \mathbf{V}\mapsto e\otimes v$ and $r_{\frak{V}}:v\in \mathbf{V}%
\mapsto v\otimes e$, indicating by $e$ the generator of $\Bbbk $.

Let us define the contravariant functor $\dagger :\mathrm{EQA}\rightarrow 
\mathrm{EQA}$, 
\begin{equation}
\dagger \,:\left( \mathbf{V},\Bbb{R}\right) \mapsto \frak{V}^{\dagger
}\doteq \left( \mathbf{V}^{\ast },\Bbb{R}^{\dagger }\right) \doteq \left( 
\mathbf{V}^{\ast },-\Bbb{R}^{\ast }\right) ,\;\;\;\;\dagger \,:l\mapsto
l^{\dagger }\doteq l^{\ast },  \label{C}
\end{equation}
being $\mathbf{V}^{\ast }$ the dual of $\mathbf{V}$, and $\Bbb{R}^{\ast }$
the transpose map w.r.t. the usual extension to $\mathbf{V}^{\otimes 2}$ of
the pairing between $\mathbf{V}$ and $\mathbf{V}^{\ast }$. It is clear that $%
\dagger ^{2}=\dagger \dagger $ is naturally equivalent to $id_{\mathrm{EQA}}$%
. In particular, $\frak{V}^{\dagger \dagger }\backsimeq \frak{V}$, $\forall 
\frak{V}\in \mathrm{EQA}$. The relation between $\boxtimes $ and $\dagger $
can be summarized by equations 
\[
\left( \frak{V}\boxtimes \frak{W}\right) ^{\dagger }\backsimeq \frak{V}%
^{\dagger }\boxtimes \frak{W}^{\dagger },\;\frak{K}^{\dagger }\backsimeq 
\frak{K}. 
\]

In terms of pairs $\left( \mathcal{A};\Bbb{R}\right) $, $\boxtimes $ and $%
\dagger $ are given by 
\[
\left( \mathcal{A};\Bbb{R}\right) \times \left( \mathcal{B};\Bbb{S}\right)
\mapsto \left( \mathcal{A}\boxtimes \mathcal{B};\Bbb{R}\boxtimes \Bbb{S}%
\right) ,\;\;\;\alpha \times \beta \mapsto \alpha \boxtimes \beta ,
\]
and 
\[
\left( \mathcal{A};\Bbb{R}\right) \mapsto \left( \mathcal{A}^{\dagger };\Bbb{%
R}^{\dagger }\right) \doteq \left( \left( \mathbf{A}_{1}^{\ast },\mathbf{A}%
_{1}^{\ast \otimes }/I\left[ \func{Im}\Bbb{R}^{\ast }\right] \right) ;-\Bbb{R%
}^{\ast }\right) ,\;\;\;\;\alpha \mapsto \alpha ^{\dagger },
\]
respectively, being $\mathcal{A}\boxtimes \mathcal{B}\doteq \left( \mathbf{A}%
_{1}\otimes \mathbf{B}_{1},\left. \left[ \mathbf{A}_{1}\otimes \mathbf{B}_{1}%
\right] ^{\otimes }\right/ \func{Im}\left[ \Bbb{R}\boxtimes \Bbb{S}\right]
\right) $. The arrows $\alpha \boxtimes \beta $ and $\alpha ^{\dagger }$ are
the extension of $\alpha _{1}\otimes \beta _{1}$ and $\alpha _{1}^{\ast }$
to an algebra map. The unit object for $\boxtimes $ is $\left( \mathcal{K};%
\Bbb{O}\right) $.

\subsection{The embedding $\mathrm{EQA}\hookrightarrow \mathrm{QA}$}

Now, we study the relationship between $\mathrm{EQA}$ and $\mathrm{QA}$ as
monoidal categories. There exists an obvious forgetful functor between these
categories.

\begin{proposition}
The function $\left( \mathcal{A};\Bbb{R}\right) \mapsto \mathcal{A}$ defines
a surjective embedding $\mathsf{U}:\mathrm{EQA}\hookrightarrow \mathrm{QA}$.
\end{proposition}

\textbf{Proof}: We just need to show the function is surjective, i.e. given $%
\mathcal{A}\in \mathrm{QA}$, there exists a compatible map $\Bbb{R}$ such
that $\mathsf{U}\left( \mathcal{A};\Bbb{R}\right) =\mathcal{A}$. Let $I\left[
\mathbf{R}\right] $ be the ideal related to $\mathcal{A}$. Consider a
decomposition $\mathbf{A}_{1}^{\otimes 2}=\mathbf{R}\oplus \mathbf{R}^{c}$,
with associated projections $\Bbb{P}$ such that $\func{Im}\Bbb{P}=\mathbf{R}$%
. Since $I\left[ \mathbf{R}\right] =I\left[ \func{Im}\Bbb{P}\right] $, then $%
\mathcal{A}\Vdash \Bbb{P}$ and the proposition follows.\ \ \ \ $\blacksquare 
$

\bigskip

On pairs $\left( \mathbf{V},\Bbb{R}\right) $ the embedding is given by $%
\left( \mathbf{V},\Bbb{R}\right) \mapsto \left( \mathbf{V},\left. \mathbf{V}%
^{\otimes }\right/ I\left[ \func{Im}\Bbb{R}\right] \right) $ (see second
part of Eq. $\left( \ref{fg}\right) $). The surjectivity is up to
isomorphisms in $\mathrm{QA}$. The functor $\mathsf{U}$ obviously preserves
the unit objects, in fact 
\[
\mathsf{U}\frak{K}=\mathsf{U}\left( \Bbbk ,\Bbb{O}\right) =\left( \Bbbk
,\left. \Bbbk ^{\otimes }\right/ I\left[ \func{Im}\Bbb{O}\right] \right)
=\left( \Bbbk ,\Bbbk ^{\otimes }\right) =\mathcal{K}, 
\]
but is not monoidal. Nevertheless,

\begin{proposition}
There exist functorial epimorphisms $\mathsf{U}\left( \frak{V}\boxtimes 
\frak{W}\right) \twoheadrightarrow \mathsf{U}\frak{V}\circ \mathsf{U}\frak{W}
$, $\frak{V},\frak{W}\in \mathrm{EQA}$, defining a natural transformation $%
\mathsf{U\,}\boxtimes \mathsf{\,}\rightarrow \circ \mathsf{\,}\left( \mathsf{%
U}\times \mathsf{U}\right) $.
\end{proposition}

\textbf{Proof:} It is clear that, given pairs $\frak{V}$ and $\frak{W}$, we
have 
\[
\begin{array}{l}
\func{Im}\Bbb{R}\boxtimes \Bbb{S}\thickapprox \varphi _{\mathbf{V},\mathbf{W}%
}\left( \func{Im}\Bbb{R}\boxtimes \Bbb{S}\right) = \\ 
\\ 
=\func{Im}\left[ \Bbb{R}\otimes \Bbb{I}+\Bbb{I}\otimes \Bbb{S}\right]
\subset \func{Im}\Bbb{R}\otimes \mathbf{W}^{\otimes 2}+\mathbf{V}^{\otimes
2}\otimes \func{Im}\Bbb{S},
\end{array}
\]
and accordingly, 
\[
\varphi _{\mathbf{V},\mathbf{W}}\left( I\left[ \func{Im}\Bbb{R}\boxtimes 
\Bbb{S}\right] \right) \subset \left( I\left[ \func{Im}\Bbb{R}\right]
\otimes \mathbf{W}^{\otimes }+\mathbf{V}^{\otimes }\otimes I\left[ \func{Im}%
\Bbb{S}\right] \right) \cap \mathbf{V}\circ \mathbf{W}. 
\]
The first ideal in above inclusion is related to the quantum space $\mathsf{U%
}\left( \frak{V}\boxtimes \frak{W}\right) $, and the latter to $\mathsf{U}%
\frak{V}\circ \mathsf{U}\frak{W}$, since it defines the algebra 
\[
\left( \left. \mathbf{V}^{\otimes }\right/ I\left[ \func{Im}\Bbb{R}\right]
\right) \circ \left( \left. \mathbf{W}^{\otimes }\right/ I\left[ \func{Im}%
\Bbb{S}\right] \right) . 
\]
Note the corresponding algebras are quotient of $\left[ \mathbf{V}\otimes 
\mathbf{W}\right] ^{\otimes }$. Hence, for every couple $\frak{V},\frak{W}%
\in \mathrm{EQA}$, we have an epimorphism $p_{\frak{V},\frak{W}}:\mathsf{U}%
\left( \frak{V}\boxtimes \frak{W}\right) \twoheadrightarrow \mathsf{U}\frak{V%
}\circ \mathsf{U}\frak{W}$.

By straightforward calculations, it can be checked commutativity of diagrams 
\begin{equation}
\begin{diagram} \QTR{sf}{U}\left( \QTR{frak}{V}\boxtimes
\QTR{frak}{W}\right) & \rOnto^{p_{\QTR{frak}{V},\QTR{frak}{W}}} &
\QTR{sf}{U}\QTR{frak}{V}\circ \QTR{sf}{U}\QTR{frak}{W} \\
\dTo^{\QTR{sf}{U}\left( \alpha \boxtimes \beta \right) }& &
\dTo_{\QTR{sf}{U}\alpha \circ \QTR{sf}{U}\beta } \\ \QTR{sf}{U}\left(
\QTR{frak}{X}\boxtimes \QTR{frak}{Y}\right) &
\rOnto^{p_{\QTR{frak}{X},\QTR{frak}{Y}}} & \QTR{sf}{U}\QTR{frak}{X}\circ
\QTR{sf}{U}\QTR{frak}{Y} \\ \end{diagram}  \label{ua}
\end{equation}
for every couple of arrows $\alpha :\frak{V}\rightarrow \frak{X}$ and $\beta
:\frak{W}\rightarrow \frak{Y}$ in $\mathrm{EQA}$.\ \ \ $\blacksquare $

\bigskip

This result, together with the ones relating monoids $\circ $ and $\otimes $%
, will be useful in order to construct the rectangular comultiplications
maps.

\section{Rectangular quantum matrix algebras}

Now, the central result. We shall show the following theorem later, in a
more general context.

\begin{theorem}
The monoidal category $\left( \mathrm{EQA},\boxtimes ,\frak{K}\right) $ is 
\textbf{rigid}, and has $\dagger $ as duality functor. For every $\frak{V}%
=\left( \mathbf{V},\Bbb{R}\right) $ in $\mathrm{EQA}$, the evaluation and
coevaluation arrows, $ev_{\frak{V}}:\frak{V}^{\dagger }\boxtimes \frak{V}%
\rightarrow \frak{K}$ and $coev_{\frak{V}}:\frak{K}\rightarrow \frak{V}%
\boxtimes \frak{V}^{\dagger }$, respectively, are given by the corresponding
maps for $\mathbf{V}$ in the rigid monoidal category $\left( \mathrm{Vct}%
_{f},\otimes ,\Bbbk \right) $. \ \ \ $\blacksquare $
\end{theorem}

We can define the internal coHom object related to a couple $\frak{W},\frak{V%
}\in \mathrm{EQA}$ as $\underline{hom}\left[ \frak{W},\frak{V}\right] \doteq 
\frak{W}^{\dagger }\boxtimes \frak{V}$, and take 
\[
\delta _{\frak{V},\frak{W}}\doteq \tau _{\frak{W},\underline{hom}\left[ 
\frak{W},\frak{V}\right] }\,\left( coev_{\frak{W}}\boxtimes I\right) \,\ell
_{\frak{V}}:\frak{V}\rightarrow \underline{hom}\left[ \frak{W},\frak{V}%
\right] \boxtimes \frak{W} 
\]
as the (left) coevaluation arrow. Its well known universality property says:
given $\frak{H}\in \mathrm{EQA}$ and $\varphi :\frak{V}\rightarrow \frak{H}%
\boxtimes \frak{W}$, there exists a unique morphism $\alpha :\underline{hom}%
\left[ \frak{W},\frak{V}\right] \rightarrow \frak{H}$ making commutative the
diagram 
\begin{equation}
\begin{diagram} & & \QTR{frak}{V} & & \\ & \ldTo^{\delta
_{\QTR{frak}{V},\QTR{frak}{W}}} & & \rdTo^{\varphi } & \\
\underline{hom}\left[ \QTR{frak}{W},\QTR{frak}{V}\right] \boxtimes
\QTR{frak}{W}& & \rTo^{\alpha \boxtimes I}&
&\QTR{frak}{H}\boxtimes\QTR{frak}{W} \\ \end{diagram}  \label{diaun}
\end{equation}
From $\left( \ref{diaun}\right) $ and general properties of monoidal
categories follow the existence of comultiplication 
\begin{equation}
\underline{hom}\left[ \frak{W},\frak{V}\right] \rightarrow \underline{hom}%
\left[ \frak{U},\frak{V}\right] \boxtimes \underline{hom}\left[ \frak{W},%
\frak{U}\right] ,\;\forall \frak{U},\frak{V},\frak{W}\in \mathrm{EQA},
\label{com}
\end{equation}
given by $\Delta _{\frak{U},\frak{V},\frak{W}}=\tau _{\underline{hom}\left[ 
\frak{W},\frak{U}\right] ,\underline{hom}\left[ \frak{U},\frak{V}\right]
}\,\left( I\boxtimes coev_{\frak{U}}\boxtimes I\right) \,\left( I\boxtimes
\ell _{\frak{V}}\right) $, and counit arrows 
\begin{equation}
\varepsilon _{\frak{V}}=ev_{\frak{V}}:\underline{end}\left[ \frak{V}\right]
\rightarrow \frak{K},\;\forall \frak{V}\in \mathrm{EQA}.  \label{coun}
\end{equation}
Coevaluations are particular comultiplications. Indeed, since $\underline{hom%
}\left[ \frak{K},\frak{V}\right] =\frak{V}$, $\forall \frak{V}\in \mathrm{EQA%
}$, it can be seen that $\delta _{\frak{V},\frak{W}}=\Delta _{\frak{W},\frak{%
V},\frak{K}}$. On the other hand, if $\frak{U}=\frak{V}=\frak{W}$, $\Delta _{%
\frak{V}}=\Delta _{\frak{V},\frak{V},\frak{V}}$ and $\varepsilon _{\frak{V}}$
gives $\underline{end}\left[ \frak{V}\right] $ a coalgebra structure in $%
\mathrm{EQA}$, and $\delta _{\frak{V}}=\delta _{\frak{V},\frak{V}}$ makes $%
\frak{V}$ an $\underline{end}\left[ \frak{V}\right] $-corepresentation in
the same category.

We add that arrows $\left( \ref{com}\right) $ and $\left( \ref{coun}\right) $
satisfy usual associativity and unit constraints, expressed by commutativity
of the following diagrams\footnote{%
Of course, we have a similar commutative diagram where $\frak{K}$ is on the
right.} 
\begin{equation}
\begin{diagram} \underline{hom}\left[ \QTR{frak}{X},\QTR{frak}{V}\right]
\boxtimes \underline{hom}\left[ \QTR{frak}{W},\QTR{frak}{X}\right] &\rTo&
\underline{hom}\left[ \QTR{frak}{Y},\QTR{frak}{V}\right] \boxtimes
\underline{hom}\left[ \QTR{frak}{X},\QTR{frak}{Y}\right] \boxtimes
\underline{hom}\left[ \QTR{frak}{W},\QTR{frak}{X}\right]\\ \uTo& &\uTo \\
\underline{hom}\left[ \QTR{frak}{W},\QTR{frak}{V}\right] &\rTo&
\underline{hom}\left[ \QTR{frak}{Y},\QTR{frak}{V}\right] \boxtimes
\underline{hom}\left[ \QTR{frak}{W},\QTR{frak}{Y}\right]\\ \end{diagram}
\label{evco}
\end{equation}
\medskip 
\begin{equation}
\begin{diagram} \QTR{frak}{K}\boxtimes \underline{hom}\left[
\QTR{frak}{W},\QTR{frak}{V}\right] & &\lTo & & \underline{end}\left[
\QTR{frak}{V}\right] \boxtimes \underline{hom}\left[
\QTR{frak}{W},\QTR{frak}{V}\right] \\ & \luTo~{\backsimeq } & &\ruTo & \\ &
& \underline{hom}\left[ \QTR{frak}{W},\QTR{frak}{V}\right] & &\\
\end{diagram}  \label{idco}
\end{equation}

\bigskip

Consider pairs $\frak{V}=\left( \mathbf{V},\Bbb{R}\right) $ and $\frak{W}%
=\left( \mathbf{W},\Bbb{S}\right) $ in $\mathrm{EQA}$, and take basis $%
\left\{ v_{i}\right\} $ and $\left\{ w_{i}\right\} $ of $\mathbf{V}$ and $%
\mathbf{W}$, respectively. The image under $\mathsf{U}$ of the internal
coHom object 
\[
\underline{hom}\left[ \frak{W},\frak{V}\right] =\frak{W}^{\dagger }\boxtimes 
\frak{V}=\left( \mathbf{W}^{\ast }\otimes \mathbf{V},\Bbb{S}^{\dagger
}\boxtimes \Bbb{R}\right) 
\]
is a quadratic algebra generated by $\mathbf{W}^{\ast }\otimes \mathbf{V}$
and obeying relations $I\left[ \func{Im}\Bbb{S}^{\dagger }\boxtimes \Bbb{R}%
\right] $. Writing $t_{i}^{j}=w^{j}\otimes v_{i}$, straightforwardly 
\begin{equation}
\func{Im}\Bbb{S}^{\dagger }\boxtimes \Bbb{R}=span\left[ \Bbb{R}%
_{ij}^{kl}\;t_{k}^{n}\otimes t_{l}^{m}-t_{i}^{k}\otimes t_{j}^{l}\;\Bbb{S}%
_{kl}^{nm}\right] _{i,j,n,m}.  \label{raa}
\end{equation}
Comparing $\left( \ref{rqa}\right) $ with $\left( \ref{raa}\right) $, we
have $\mathsf{U}\underline{hom}\left[ \frak{W},\frak{V}\right] =A\left( \Bbb{%
R}:\Bbb{S}\right) $. In particular, $\mathsf{U}\underline{end}\left[ \frak{V}%
\right] =A\left( \Bbb{R}\right) $. Thus, the algebras $A\left( \Bbb{R}:\Bbb{S%
}\right) $ (resp. bialgebras $A\left( \Bbb{R}\right) $) are the coordinate
ring of an equipped quantum space with structure $\Bbb{S}^{\dagger
}\boxtimes \Bbb{R}$ (resp. $\Bbb{R}^{\dagger }\boxtimes \Bbb{R}$),
representing the space of homomorphisms from $\frak{W}^{op}$ to $\frak{V}%
^{op}$.

\subsection{Rectangular comultiplication and counit maps}

Now, we are going to construct rectangular comultiplication and counit maps
defined in \cite{mm} for algebras $A\left( \Bbb{R}:\Bbb{S}\right) $. This
can be done in steps below:

\begin{enumerate}
\item  Apply the functor $\mathsf{U}$ to comultiplications given in $\left( 
\ref{com}\right) $ to obtain the maps 
\[
\mathsf{U}\underline{hom}\left[ \frak{W},\frak{V}\right] \rightarrow \mathsf{%
U}\left( \underline{hom}\left[ \frak{U},\frak{V}\right] \boxtimes \underline{%
hom}\left[ \frak{W},\frak{U}\right] \right) .
\]

\item  Compose it with the functorial epimorphism 
\[
\mathsf{U}\left( \underline{hom}\left[ \frak{U},\frak{V}\right] \boxtimes 
\underline{hom}\left[ \frak{W},\frak{U}\right] \right) \twoheadrightarrow 
\mathsf{U}\underline{hom}\left[ \frak{U},\frak{V}\right] \circ \mathsf{U}%
\underline{hom}\left[ \frak{W},\frak{U}\right] .
\]

\item  Apply the forgetful functor $\mathsf{F}:\left( \mathbf{A}_{1},\mathbf{%
A}\right) \mapsto \mathbf{A}$ to this composition. This gives us an algebra
homomorphism 
\[
\mathsf{FU}\underline{hom}\left[ \frak{W},\frak{V}\right] \rightarrow 
\mathsf{F}\left( \mathsf{U}\underline{hom}\left[ \frak{U},\frak{V}\right]
\circ \mathsf{U}\underline{hom}\left[ \frak{W},\frak{U}\right] \right) .
\]

\item  Finally, compose the latter with the functorial inclusion 
\[
\mathsf{F}\left( \mathsf{U}\underline{hom}\left[ \frak{U},\frak{V}\right]
\circ \mathsf{U}\underline{hom}\left[ \frak{W},\frak{U}\right] \right)
\hookrightarrow \mathsf{FU}\underline{hom}\left[ \frak{U},\frak{V}\right]
\otimes \mathsf{FU}\underline{hom}\left[ \frak{W},\frak{U}\right] .
\]
\end{enumerate}

The resulting maps are precisely the arrows $A\left( \Bbb{R}:\Bbb{S}\right)
\rightarrow A\left( \Bbb{R}:\Bbb{T}\right) \otimes A\left( \Bbb{T}:\Bbb{S}%
\right) $ defined in \cite{mm}. For counits:

\begin{enumerate}
\item  Apply $\mathsf{FU}$ to the counit arrow $\varepsilon _{\frak{V}}:%
\underline{end}\left[ \frak{V}\right] \rightarrow \frak{K}$ given in $\left( 
\ref{coun}\right) $, to obtain the arrow $A\left( \Bbb{R}\right) \rightarrow
\Bbbk ^{\otimes }$.

\item  Compose to surjection $\Bbbk ^{\otimes }\twoheadrightarrow \Bbbk $.
\end{enumerate}

The algebra homomorphism $A\left( \Bbb{R}\right) \rightarrow \Bbbk $ we
obtain, together with comultiplication $A\left( \Bbb{R}\right) \rightarrow
A\left( \Bbb{R}\right) \otimes A\left( \Bbb{R}\right) $ gives $A\left( \Bbb{R%
}\right) $ a bialgebra structure. In fact, diagrams $\left( \ref{evco}%
\right) $ and $\left( \ref{idco}\right) $, properly combined with $\left( 
\ref{fa}\right) $, $\left( \ref{fu}\right) $ and $\left( \ref{ua}\right) $,
lead us to associativity and unit constraints for arrows 
\[
A\left( \Bbb{R}:\Bbb{S}\right) \rightarrow A\left( \Bbb{R}:\Bbb{T}\right)
\otimes A\left( \Bbb{T}:\Bbb{S}\right) ,\;\;A\left( \Bbb{R}\right)
\rightarrow \Bbbk . 
\]

Summing up, we have constructed square and rectangular quantum matrices as
internal coHom objects in the rigid category of equipped quantum spaces,
giving a generalization of \emph{FRT} construction in the scenario of Manin
quantum groups.

\subsection{The inclusions $\mathsf{U}\protect\underline{hom}\left[ \frak{W},%
\frak{V}\right] ^{op}\hookrightarrow \protect\underline{hom}\left[ \mathsf{U}%
\frak{W},\mathsf{U}\frak{V}\right] ^{op}$}

Let us show there exist epis $\underline{hom}\left[ \mathsf{U}\frak{W},%
\mathsf{U}\frak{V}\right] \twoheadrightarrow A\left( \Bbb{R}:\Bbb{S}\right) $%
, with $\underline{hom}\left[ \mathsf{U}\frak{W},\mathsf{U}\frak{V}\right] $
the internal coHom object in $\mathrm{QA}$ related to quantum spaces $%
\mathsf{U}\frak{W}$ and $\mathsf{U}\frak{V}$. As we have claimed in the
introduction, this is an indication that quadratic algebras $A\left( \Bbb{R}:%
\Bbb{S}\right) $ represent homomorphisms between structurally richer spaces,
w.r.t. coHom objects in $\mathrm{QA}$.

Consider the evaluation map $\delta _{\frak{V},\frak{W}}:\frak{V}\rightarrow 
\underline{hom}\left[ \frak{W},\frak{V}\right] \boxtimes \frak{W}$ in $%
\mathrm{EQA}$. In the previously given basis $\left\{ v_{i}\right\} $ and $%
\left\{ w_{i}\right\} $ this map is defined by the assignment $v_{i}\mapsto
t_{i}^{j}\otimes w_{j}$, putting again $t_{i}^{j}=w^{j}\otimes v_{i}$. It is
functored by $\mathsf{U}$ to an arrow $\mathsf{U}\frak{V}\rightarrow \mathsf{%
U}\left( \underline{hom}\left[ \frak{W},\frak{V}\right] \boxtimes \frak{W}%
\right) $ which, composed to the functorial epi 
\[
\mathsf{U}\left( \underline{hom}\left[ \frak{W},\frak{V}\right] \boxtimes 
\frak{W}\right) \twoheadrightarrow \mathsf{U}\underline{hom}\left[ \frak{W},%
\frak{V}\right] \circ \mathsf{U}\frak{W}, 
\]
defines in $\mathrm{QA}$ another arrow $\varphi :\mathsf{U}\frak{V}%
\rightarrow \mathsf{U}\underline{hom}\left[ \frak{W},\frak{V}\right] \circ 
\mathsf{U}\frak{W}$, also defined by $v_{i}\mapsto t_{i}^{j}\otimes w_{j}$.
From universality of internal coHom objects in $\mathrm{QA}$ (see Equation $%
\left( \ref{diaun}\right) $), there exists a unique arrow 
\[
\alpha :\underline{hom}\left[ \mathsf{U}\frak{W},\mathsf{U}\frak{V}\right]
\rightarrow \mathsf{U}\underline{hom}\left[ \frak{W},\frak{V}\right] 
\]
such that $\varphi =\left( \alpha \circ I\right) \,\delta _{\mathsf{U}\frak{V%
},\mathsf{U}\frak{W}}$, being $\delta _{\mathsf{U}\frak{V},\mathsf{U}\frak{W}%
}$ the coevaluation associated to $\underline{hom}\left[ \mathsf{U}\frak{W},%
\mathsf{U}\frak{V}\right] \in \mathrm{QA}$. Recalling that (since $\mathsf{U}%
\frak{W}$ and $\mathsf{U}\frak{V}$ are quadratic algebras generated by $%
\mathbf{W}$ and $\mathbf{V}$, resp.) $\underline{hom}\left[ \mathsf{U}\frak{W%
},\mathsf{U}\frak{V}\right] $ is generated by $\mathbf{W}^{\ast }\otimes 
\mathbf{V}$ and $\delta _{\mathsf{U}\frak{V},\mathsf{U}\frak{W}}\left(
v_{i}\right) =t_{i}^{j}\otimes w_{j}$, then $\alpha $ is the identity on
generators. Consequently, due to $\mathsf{U}\underline{hom}\left[ \frak{W},%
\frak{V}\right] $ is also generated by $\mathbf{W}^{\ast }\otimes \mathbf{V}$%
, $\alpha $ is an algebra epimorphism.

Hence, we have shown: For every couple $\frak{W},\frak{V}$ of equipped
quantum space, the quantum space $\underline{hom}\left[ \mathsf{U}\frak{W},%
\mathsf{U}\frak{V}\right] ^{op}$ `contains' $\mathsf{U}\underline{hom}\left[ 
\frak{W},\frak{V}\right] ^{op}$ as a subspace.

\section{The equipped conic quantum spaces}

Given a finite dimensional $\Bbbk $-vector space\textbf{\ }$\mathbf{V}\in $%
\textrm{Vct}$_{f}$, consider the degree cero homogeneous linear
endomorphisms of $\mathbf{V}^{\otimes }$, i.e. 
\[
\Bbb{R}\in End_{\mathrm{Vct}}\left[ \mathbf{V}^{\otimes }\right]
\;\;such\;that\;\;\Bbb{R}\left( \mathbf{V}^{\otimes n}\right) \subset 
\mathbf{V}^{\otimes n}.
\]
Of course, each map $\Bbb{R}$ is defined by a family $\left\{ \Bbb{R}%
_{n}\right\} _{n\in \Bbb{N}_{0}}$ of linear maps $\Bbb{R}_{n}:\mathbf{V}%
^{\otimes n}\rightarrow \mathbf{V}^{\otimes n}$. In terms of these
endomorphisms, all above constructions can be repeated word by word in the
category $\mathrm{CA}$. That is to say, we can define \emph{equipped conic
quantum spaces }as pairs $\left( \mathcal{A},\Bbb{R}\right) $, $\mathcal{A}%
\in \mathrm{CA}$, such that $\mathcal{A}\Vdash \Bbb{R}$, i.e. $\ker \left[ 
\mathbf{A}_{1}^{\otimes }\twoheadrightarrow \mathbf{A}\right] =I\left[ \func{%
Im}\Bbb{R}\right] $. We just must change the defining condition for
morphisms $\left( \mathcal{A},\Bbb{R}\right) \rightarrow \left( \mathcal{B},%
\Bbb{S}\right) $ by $\alpha _{1}^{\otimes }\Bbb{\,R}=\Bbb{S\,}\alpha
_{1}^{\otimes }$, or $\alpha _{1}^{\otimes n}\Bbb{\,R}_{n}=\Bbb{S}_{n}\Bbb{\,%
}\alpha _{1}^{\otimes n}$ for all $n\in \Bbb{N}_{0}$. Let us call $\mathrm{%
ECA}$ the related category. It can be shown, the category of pairs $\frak{V}%
=\left( \mathbf{V},\Bbb{R}\right) $ such that $\func{Im}\Bbb{R}_{1}=\left\{
0\right\} $, with morphisms given by linear maps $l$ such that $l^{\otimes }%
\Bbb{\,R}=\Bbb{S\,}l^{\otimes }$, is equivalent to $\mathrm{ECA}$. The
functors $\boxtimes $ and $\dagger $ are defined by Eqs. $\left( \ref{A}%
\right) $, $\left( \ref{B}\right) $ and $\left( \ref{C}\right) $. $\mathrm{%
EQA}$ can be seen as a full subcategory of $\mathrm{ECA}$ by regarding the
endomorphisms $\Bbb{R}:\mathbf{V}^{\otimes 2}\rightarrow \mathbf{V}^{\otimes
2}$ as the map $\Bbb{R}_{2}$ of an homogeneous endomorphism $\Bbb{R}\in End_{%
\mathrm{Vct}}\left[ \mathbf{V}^{\otimes }\right] $ such that $\Bbb{R}_{n}$
is the null map if $n\neq 2$. This enable us to define the full subcategory
of equipped $m$-th quantum spaces $\mathrm{ECA}^{m}$ in terms of pairs $%
\left( \mathbf{V},\Bbb{R}\right) $ with $\Bbb{R}_{n}$ the null map for all $%
n\neq m$.

Again, function $\left( \mathbf{V},\Bbb{R}\right) \mapsto \left( \mathbf{V}%
,\left. \mathbf{V}^{\otimes }\right/ I\left[ \func{Im}\Bbb{R}\right] \right) 
$ gives rise to a surjective embedding $\mathsf{U}$, now in $\mathrm{CA}$,
such that the functorial epimorphisms $\mathsf{U}\left( \frak{V}\boxtimes 
\frak{W}\right) \twoheadrightarrow \mathsf{U}\frak{V}\circ \mathsf{U}\frak{W}
$ are still valid.

Now, let us prove \textbf{Theor. 4} in this more general setting.

\begin{theorem}
The category of equipped conic quantum spaces is \textbf{rigid }w.r.t. the
monoidal structure $\boxtimes $, and has $\dagger $ as duality functor.
\end{theorem}

\textbf{Proof:}\ Consider an object $\frak{V}=\left( \mathbf{V},\Bbb{R}%
\right) $. We define the evaluation and the coevaluation morphisms, 
\[
ev_{\frak{V}}:\frak{V}^{\dagger }\boxtimes \frak{V}\rightarrow \frak{K}%
\;\;\;\;and\;\;\;\;coev_{\frak{V}}:\frak{K}\rightarrow \frak{V}\boxtimes 
\frak{V}^{\dagger },
\]
as the usual pairing $v\otimes v^{\prime }\mapsto \left\langle v,v^{\prime
}\right\rangle \,e$ and coevaluation of $\mathbf{V}$ and $\mathbf{V}^{\ast }$%
, respectively, being $e$ the generator of $\Bbbk $. We must show that they
are effectively arrows in $\mathrm{E}$\textrm{$CA$}.

Note that the tensor product map $ev_{\frak{V}}^{\otimes }$ defines the
algebra homomorphism $\mathbf{V}^{\ast }\circ \mathbf{V}\rightarrow \Bbbk
^{\otimes }$, 
\[
u\otimes v\in \mathbf{V}^{\ast \otimes n}\otimes \mathbf{V}^{\otimes
n}\mapsto \left\langle u,v\right\rangle \;e^{n}\in \Bbbk ^{\otimes n}, 
\]
where we use $\varphi _{\mathbf{V}^{\ast },\mathbf{V}}$ to identify the
algebras $\mathbf{V}^{\ast }\circ \mathbf{V}$ and $\left[ \mathbf{V}^{\ast
}\otimes \mathbf{V}\right] ^{\otimes }$. By direct calculations (and from
the very definition of $\Bbb{R}^{\dagger }$), it can be seen 
\[
ev_{\frak{V}}^{\otimes }\,\,\Bbb{R}^{\dagger }\boxtimes \Bbb{R}=0=\Bbb{O}%
\,ev_{\frak{V}}^{\otimes }. 
\]
To show the analogous equation for the coevaluation map, let us first
introduce some notation. Let $\left\{ v_{i}\right\} $ be a basis of $\mathbf{%
V}$. Construct for each $n\in \Bbb{N}$ a basis $\left\{ v_{R}\right\} $ of $%
\mathbf{V}^{\otimes n}$, being $R=\left( r_{1},...,r_{n}\right) $ a
multi-index with $1\leq r_{k}\leq \dim \mathbf{V}$, $\forall k=1,...,n$, in
such a way that $v_{R}=v_{r_{1}}\otimes ...\otimes v_{r_{n}}.$ Consider also
the basis $\left\{ v^{R}\right\} $ of $\mathbf{V}^{\ast \otimes n}$ dual to $%
\left\{ v_{R}\right\} $ w.r.t. the usual pairing $\left\langle \cdot ,\cdot
\right\rangle _{n}:\mathbf{V}^{\ast \otimes n}\otimes \mathbf{V}^{\otimes
n}\rightarrow \Bbbk $. In these terms, the algebra homomorphism 
\[
e^{n}\in \Bbbk ^{\otimes }\mapsto v_{R}\otimes v^{R}\in \mathbf{V}\circ 
\mathbf{V}^{\ast } 
\]
(sum over repeated (multi)indices is understood) coincides with the map $%
coev_{\frak{V}}^{\otimes }$. To see that equation 
\[
\Bbb{R}\boxtimes \Bbb{R}^{\dagger }\,coev_{\frak{V}}^{\otimes }=coev_{\frak{V%
}}^{\otimes }\,\Bbb{O} 
\]
holds, we just have to prove $\Bbb{R}\otimes I-I\otimes \Bbb{R}^{\ast }$
evaluated on $v_{R}\otimes v^{R}$ is equal to cero, i.e. 
\[
\left[ \Bbb{R}\left( v_{R}\right) \otimes v^{R}-v_{S}\otimes \Bbb{R}^{\ast
}\left( v^{S}\right) \right] =0. 
\]
Writing $\Bbb{R}\left( v_{R}\right) =\Bbb{R}_{R}^{S}\;v_{S}$, we have for
the transpose $\Bbb{R}^{\ast }\left( v^{S}\right) =v^{R}\;\Bbb{R}_{R}^{S}$,
and consequently the left member of equation above is identically cero. So,
the theorem is proven.\ \ \ $\blacksquare $

\subsubsection*{Acknowledgments}

The author thanks to CNEA and Fundaci\'{o}n Antorchas, Argentina, for
financial support.

\end{document}